\documentclass[12pt]{article}
\usepackage{amsfonts,amssymb,amsmath}
\newcommand{\beq}{ \begin{equation}}
\newcommand{\eeq}{ \end{equation}}

\begin{document}
\pagestyle{empty}

\title{Markov Stochastic Operators of Heredity}

\vskip 1cm \maketitle
\begin{center}
N. N. Ganikhodjaev \\
\end{center}
\begin{center}

Kulliyyah of Science, International Islamic University Malaysia,
\\25200 Kuantan, Malaysia and \\
Institute of Mathematics,Tashkent,700125,Uzbekistan.
\end{center}
\noindent

\section {Quadratic Stochastic Operator}
   Let
$$ S^{m-1}=\{\textbf{x}=(x_1,\cdots,x_m)\in R^m: x_i\geq 0, \sum^m_{i=1}x_i=1 \} $$
is the $(m-1)$ -dimensional canonical simplex in $R^m.$ The
transformation $V:S^{m-1}\rightarrow S^{m-1}$ is called a quadratic
stochastic operator(q.s.o.) if
$$ V: (V\textbf{x})_k=\sum^m_{i,j=1}p_{ij,k}x_ix_j, \quad  (k=1,\cdots,m)\eqno(1) $$
where
\[
 a)p_{ij,k}\geq0;\  \quad
  b)p_{ij,k}=p_{ji,k};\ \quad\mbox {and}\quad
  c)\sum^m_{k=1}p_{ij,k}=1\,
\]
 for arbitrary $ i,j,k=1,\cdots,m.$
 Quadratic stochastic operator was first introduced in [1].Such
operator frequently arises in many models of mathematical
genetics[1-3].Consider a biological population,that is a community
of organisms closed with respect to reproduction.Assume that each
individual in this population belongs to precisely one species
 $1,\cdots, m . $ The scale of species is such that the species of the
parents $ i $ and $ j $ unambiguously determines the probability
of every species $k$ for the first generation of direct
descendants.Denote this probability,that is to be called the
heredity coefficient, by $p_{ij,k}.$ It is then obvious that
$p_{ij,k}\geq 0$ for all $i,j,k$ and that
$$\sum^m_{k=1}p_{ij,k}=1  (i,j,k=1,\cdots,m). $$

Assume that the population is so large that frequency fluctuations
can be neglected.Then the state of the population can be described
by the tuple $(x_1,x_2,\cdots,x_m)$ of species probabilities,that is
$ x_k$ is the fraction of the species $k$ in the total population.In
the case of panmixia (random interbreeding) the parent pairs $i$ and
$j$ arise for a fixed state $\textbf{x}=(x_1,x_2,\cdots,x_m)$ with
probability $x_ix_j.$Hence the total probability of the species $k$
in the first generation of direct descendants is defined by q.s.o.
(1).

 \section {Markov Stochastic Operator}
 In the case of panmixia (random interbreeding) the parent pairs $i$ and
$j$ arise for a fixed state $\textbf{x}=(x_1,x_2,\cdots,x_m)$ with
probability $x_ix_j.$ In this case we call corresponding
transformation (1)
$$ V: (V\textbf{x})_k=\sum^m_{i,j=1}p_{ij,k}x_ix_j, \quad  (k=1,\cdots,m) $$
a \emph{Bernoulli q.s.o. }

Let $\Pi=(q_{ij})_{i,j=1}^m $ be a stochastic matrix and
$\textbf{x}=(x_1,x_2,\cdots,x_m)$ be a fixed state of
population.If parents pairs $i$ and $j$ arise with probability
$x_iq_{ij}$  we  call such interbreeding  $\Pi$-\emph{panmixia} or
$\Pi$-\emph{random interbreeding}. Under $\Pi$-panmixia the total
probability of the species $k$ in the first generation of direct
descendants is defined as

$$ V_{\Pi}: (V_{\Pi}\textbf{x})_k=\sum^m_{i,j=1}p_{ij,k}q_{ij}x_i \ \ (k=1,\cdots,m).\eqno(2) $$

 We call a
transformation (2) $ V_{\Pi}:S^{m-1}\rightarrow S^{m-1}$  a
\emph{Markov quadratic stochastic operator of heredity.}

 Under $\Pi$-panmixia the total probability of
the species $ k $ in the $ (n+1)th $ generation of direct
descendants is defined as
$$ x^{(n+1)}_k=(V^n_{\Pi}(\textbf{x}))_k= \sum^m_{i,j=1}p^{(n,n+1)}_{ij,k}q^{(n,n+1)}_{ij}x^{(n)}_i \ \
(k=1,\cdots,m), $$ where $ n=0,1,2,\cdots,$ and
$\textbf{x}^{(0)}=\textbf{x}.$

Here we have two type of nonhomogeneity :nonhomogeneity with respect
to stochastic  cubic matrices $P^{(n,n+1)}=\parallel
p^{(n,n+1)}_{ij,k}\parallel$ and nonhomogeneity with respect to
stochastic matrices of random interbreeding $\Pi^{(n,n+1)}=\parallel
q^{(n,n+1)}_{ij}\parallel_{i,j=1}^m. $

 If  $\Pi^{(n,n+1)}=\parallel q^{(n,n+1)}_{ij}\parallel_{i,j=1}^m $ is defined as $q^{(n,n+1)}_{ij}=x^{(n)}_j$
then Markov quadratic stochastic operator of heredity(2) is reduced
to Bernoulli q.s.o. (1).
\section{Non-homogeneous Markov Chains}
Let $SM_m$ be a set of stochastic matrices $\Pi=(q_{ij})_{i,j=1}^m
.$ Given cubic matrix $P =\parallel p_{ij,k}\parallel$ we can write
in the following  form $P=(P_1\mid P_2 \mid \cdots \mid P_m)$ where
a stochastic matrix $P_i$ has following form for any $i$
\begin{displaymath}
P_i=
 \left(
 \begin{array} {cccc}
 p_{i1,1} & p_{i1,2}&\cdots & p_{i1,m}\\
 p_{i2,1} & p_{i2,2} &  \cdots & p_{i2,m} \\

 \vdots &\vdots & \vdots & \vdots  \\
 p_{im,1}&p_{im,2}&\cdots &p_{im,m}
 \end {array}\right )
 \end{displaymath}

 Let
 \begin{displaymath}
 \Pi=
 \left(
 \begin{array} {cccc}
 q_{11} & q_{12} &\cdots & q_{1,m}\\
 q_{21} & q_{22} & \cdots & q_{2m} \\

 \vdots &\vdots & \vdots & \vdots  \\
 q_{m1}& q_{m2}& \cdots &q_{m,m}
 \end {array}\right )
 \end{displaymath}
be a stochastic matrix of random interbreeding that belong to
$SM_m$.Assume $\textbf{q}_k=(q_{k1},q_{k2},\cdots,q_{km})$.

It is evident $\textbf{q}_k \in S^{m-1}$ for arbitrary $ k $. For
given cubic matrix $P$ let us define an operator $\textbf{P}$ acting
on $SM_m$:
 \begin{displaymath}
 \Pi \textbf{P}=
 \left(
 \begin{array} {c}
 \textbf{q}_1P_1\\
 \textbf{q}_2P_2 \\

 \vdots   \\
 \textbf{q}_m P_m
 \end {array}\right )\eqno(3)
 \end{displaymath}
Here $\textbf{q}_k P_k\in S^{m-1}$ is the result of action operator
$P_k$ to vector $\textbf{q}_k.$ Evidently an operator $\textbf{P}$
transform $SM_m$ into itself. Thus for arbitrary
$\textbf{x}=(x_1,x_2,\cdots,x_m)\in S^{m-1}$ a Markov stochastic
operator of heredity (2) has following form:
$$V_{\Pi}\textbf{x}= \textbf{x}(\Pi \textbf{P}). \eqno(4)$$
According (3) $\Pi \textbf{P}$ is a stochastic matrix,so that the
trajectory of a Markov stochastic operator of heredity is defined as
a trajectory of nonhomogeneous Markov chain. Thus the study of the
asymptotic behavior of Markov stochastic operator of heredity
correspond to study of asymptotic behavior of a non-homogeneous
Markov chain $Q(1),Q(2),\cdots $,where $Q(k)$ is one -step
transition matrix with time step $k$. Note that
$$Q(k)=\Pi^{(k,k+1)}P^{(k,k+1)}.\eqno(5)$$
{\bf Remark 1} For Bernoulli q.s.o. one-step transition matrix
$Q(k)$ with time step $k$ is defined by initial distribution
$\textbf{x}^{(0)}=\textbf{x}.$ Thus for any fixed
$\textbf{x}^{(0)}=\textbf{x}$ we have one-step transition matrix
$Q(k)(\textbf{x})$ .

\section{Ergodic Theorem}
We will say that the ergodic theorem holds for transformation $
V_{\Pi}:S^{m-1}\rightarrow S^{m-1}$
  if for each initial point $ \textbf{x}\in
S^{m-1}$ the limit
$$ \lim _{k\rightarrow \infty} \frac{1}{k}\sum_{n=0}^{k-1}  V_{\Pi}^n (\textbf{x})$$
 exists.

On the basis of numerical calculations Ulam conjectured [5] that the
ergodic theorem holds for any Bernoulli q.s.o.. In 1977 Zakharevich
[5] showed that this conjecture is false in general. He proved that
the ergodic theorem does not hold for q.s.o. $V$ , which is defined
on the simplex $S^2$ by the formula
$$ V:(x,y,z)\rightarrow (x^2 +2xy,y^2 +2yz,z^2 +2xz) \eqno (6)$$

\section{Ergodic non-homogeneous Markov Chains}
Given a non-homogeneous Markov chain $(Q(1),Q(2),\cdots )$, define

$$Q^{i:j}=Q(i+1)Q(i+2)\cdots Q(j).$$

The Markov chain is said to be weak ergodic if for any given
$i\geq0$,as $j\rightarrow \infty,$ each column of $Q^{i:j}$ gets to
be a constant column.

(Though the constant may change with $j$.)Unlike in the homogeneous
case,generally one should not expect that
$\lim_{j\rightarrow\infty}Q^{i:j}$ exists.

Hajnal [7]introduced the powerful concept-scrambling matrix.

 {\bf Definition 1 }(Scrambling). Let $Q$ be the one-step transition
 matrix of a given Markov chain. $Q$ is said to be a scrambling
 matrix if for any two distinct states $i$ and $j$ , there always
 exists a state $k$, such that both one-step transitions are
 possible:$i\rightarrow k$, and $j\rightarrow k$; or equivalently,
 $q_{ik}$ and $q_{jk}$ are both positive.
 In another word, for each pair of rows in $Q$, there is a column
 such that the intersection are both positive.

Let $\Sigma$ be a set of stochastic matrices. A $\Sigma$-Markov
chain is one whose transition matrices are all
 taken from $\Sigma.$

 Jianhong Shen[6] proved following statement.

 {\bf Proposition 1} Let $\Sigma$ be a compact set of scrambling stochastic matrices.
 Then any $\Sigma$-Markov chain is weak ergodic.

  \section {Zakharevich's example}
Let us consider Bernoulli q.s.o.
$$ V:(x,y,z)\rightarrow (x^2 +2xy,y^2 +2yz,z^2 +2xz) .$$
Here corresponding cubic matrix has following form $P=(P_1\mid P_2
\mid P_3)$
\begin{displaymath}
P_1=
 \left(
 \begin{array} {ccc}
 1 & 0 & 0\\
 1 &0 & 0 \\
 0 & 0& 1
 \end {array}\right )
 \end{displaymath}

\begin{displaymath}
P_2=
 \left(
 \begin{array} {ccc}
 1 & 0 & 0\\
 0 &1 &0 \\
 0& 1 & 0
 \end {array}\right )
 \end{displaymath}

\begin{displaymath}
P_3=
 \left(
 \begin{array} {ccc}
 0 &0 & 1 \\
 0 &1 & 0 \\
 0 & 0 & 1
 \end {array}\right )
 \end{displaymath}
Then for any fixed $\textbf{x}^{(0)}=(x,y,z)$ one-step transition
matrix $Q(k)(\textbf{x})$ has following form
\begin{displaymath}
Q(k)(\textbf{x})=
 \left(
 \begin{array} {ccc}
 x^{(k)}+y^{(k)} & 0 & z^{(k)}\\
 x^{(k)} &y^{(k)}+z^{(k)} & 0 \\
 0 & y^{(k)}& x^{(k)}+z^{(k)}
 \end {array}\right )
 \end{displaymath}
where $(x^{(k+1)},y^{(k+1)},z^{(k+1)})=V(x^{(k)},y^{(k)},z^{(k)})$
for $k=0,1,\cdots $

It is evidently that for any $k$ and any $(x,y,z)\in S^2$ a
stochastic matrix $Q(k)(\textbf{x})$ is the scrambling matrix.

Using a geometric approach to ergodic non-homogeneous Markov chains
developed by Jianhong Shen we can prove the following
proposition.

{\bf Proposition 2 } The  set of all stochastic matrices $
\{Q(k)(\textbf{x}):\textbf{x}\in S^2,k=1,2,\cdots\}$ is a compact
set.

{\bf Corollary} Considered by Zakharevich Bernoulli quadratic
stochastic operator generate non-homogeneous weak ergodic Markov
chains.

Proof follows from Proposition 1.

{\bf Conclusion} For any $\textbf{x}\in S^2$ corresponding
non-homogeneous Markov chain $\{Q(k)(\textbf{x}):k=1,2,\cdots\} $ is
weak ergodic.

Thus the continual family of weak ergodic non-homogeneous Markov
chains is constructed.

 {\bf References}

1.Bernstein S.N.,Solution of a mathematical problem connected with
the theory of heredity,{\it Ann.Sci. de l'Ukraine},1:83- 114(1924)

2.Kesten H., Quadratic transformations: a model for population
growth.I, II,{\it Adv.Appl.Prob.} ,{\bf 2}, 1-82, 179-228 (1970) .

3.Lyubich Yu.I., Basic concepts and theorems of the evolution
genetics of free populations , {\it Russian Math.Surveys}, {\bf
26}:5,51-116 (1971).

4.Ulam S.,{\it A collection of mathematical problems}, Interscience
Publishers,New-York-London 1960.

5.Zakharevich M.I.,On behavior of trajectories and the ergodic
hypothesis for quadratic transformations of the simplex,{\it Russian
Math.Surveys} {\bf33}:6,265-266 (1978)

6.Jianhong Shen, A geometric approach to ergodic non-homogeneous
Markov chains,in Wavelet Analysis and Multiresolution Methods, Ed.
in Wavelet Analysis and Multiresolution Methods, Ed. T.-X. He,
Lecture Notes in Pure and Applied Mathematics, 212, Marcel Dekker,
pp. 341-366, (2000).

7.Hajnal J.,Weak ergodicity in nonhomogeneous Markov chains.{\it
Proc.Camb.Phl.Soc.}54:233-246,(1958).

\end{document}